\documentclass[10pt,a4paper,twoside]{article}
\usepackage{a4,amsfonts}

\newread\epsffilein    % file to \read
\newif\ifepsffileok    % continue looking for the bounding box?
\newif\ifepsfbbfound   % success?
\newif\ifepsfverbose   % report what you're making?
\newif\ifepsfdraft     % use draft mode?
\newdimen\epsfxsize    % horizontal size after scaling
\newdimen\epsfysize    % vertical size after scaling
\newdimen\epsftsize    % horizontal size before scaling
\newdimen\epsfrsize    % vertical size before scaling
\newdimen\epsftmp      % register for arithmetic manipulation
\newdimen\pspoints     % conversion factor
\def\epsfbox#1{\global\def\epsfllx{72}\global\def\epsflly{72}%
   \global\def\epsfurx{540}\global\def\epsfury{720}%
   \def\lbracket{[}\def\testit{#1}\ifx\testit\lbracket
   \let\next=\epsfgetlitbb\else\let\next=\epsfnormal\fi\next{#1}}%
\def\epsfgetlitbb#1#2 #3 #4 #5]#6{\epsfgrab #2 #3 #4 #5 .\\%
   \epsfsetgraph{#6}}%
\def\epsfnormal#1{\epsfgetbb{#1}\epsfsetgraph{#1}}%
\def\epsfgetbb#1{%
\openin\epsffilein=#1
\ifeof\epsffilein\errmessage{I couldn't open #1, will ignore it}\else
   {\epsffileoktrue \chardef\other=12
    \def\do##1{\catcode`##1=\other}\dospecials \catcode`\ =10
    \loop
       \read\epsffilein to \epsffileline
       \ifeof\epsffilein\epsffileokfalse\else
          \expandafter\epsfaux\epsffileline:. \\%
       \fi
   \ifepsffileok\repeat
   \ifepsfbbfound\else
    \ifepsfverbose\message{No bounding box comment in #1; using defaults}\fi\fi
   }\closein\epsffilein\fi}%
\def\epsfclipoff{\def\epsfclipstring{\ifepsfdraft\space clip\fi}}%
\epsfclipoff
\def\epsfsetgraph#1{%
   \epsfrsize=\epsfury\pspoints
   \advance\epsfrsize by-\epsflly\pspoints
   \epsftsize=\epsfurx\pspoints
   \advance\epsftsize by-\epsfllx\pspoints
   \epsfxsize\epsfsize\epsftsize\epsfrsize
   \ifnum\epsfxsize=0 \ifnum\epsfysize=0
      \epsfxsize=\epsftsize \epsfysize=\epsfrsize
      \epsfrsize=0pt
     \else\epsftmp=\epsftsize \divide\epsftmp\epsfrsize
       \epsfxsize=\epsfysize \multiply\epsfxsize\epsftmp
       \multiply\epsftmp\epsfrsize \advance\epsftsize-\epsftmp
       \epsftmp=\epsfysize
       \loop \advance\epsftsize\epsftsize \divide\epsftmp 2
       \ifnum\epsftmp>0
          \ifnum\epsftsize<\epsfrsize\else
             \advance\epsftsize-\epsfrsize \advance\epsfxsize\epsftmp \fi
       \repeat
       \epsfrsize=0pt
     \fi
   \else \ifnum\epsfysize=0
     \epsftmp=\epsfrsize \divide\epsftmp\epsftsize
     \epsfysize=\epsfxsize \multiply\epsfysize\epsftmp   
     \multiply\epsftmp\epsftsize \advance\epsfrsize-\epsftmp
     \epsftmp=\epsfxsize
     \loop \advance\epsfrsize\epsfrsize \divide\epsftmp 2
     \ifnum\epsftmp>0
        \ifnum\epsfrsize<\epsftsize\else
           \advance\epsfrsize-\epsftsize \advance\epsfysize\epsftmp \fi
     \repeat
     \epsfrsize=0pt
    \else
     \epsfrsize=\epsfysize
    \fi
   \fi
   \ifepsfverbose\message{#1: width=\the\epsfxsize, height=\the\epsfysize}\fi
   \epsftmp=10\epsfxsize \divide\epsftmp\pspoints
   \vbox to\epsfysize{\vfil\hbox to\epsfxsize{%
      \ifnum\epsfrsize=0\relax
        \includegraphics{\ifepsfdraft}%
      \else
        \epsfrsize=10\epsfysize \divide\epsfrsize\pspoints
        \includegraphics{\ifepsfdraft}%
      \fi
      \hfil}}%
\global\epsfxsize=0pt\global\epsfysize=0pt}%
{\catcode`\%=12 \global\let\epsfpercent=%\global\def\epsfbblit{%BoundingBox}}%
\long\def\epsfaux#1#2:#3\\{\ifx#1\epsfpercent
   \def\testit{#2}\ifx\testit\epsfbblit
      \epsfgrab #3 . . . \\%
      \epsffileokfalse
      \global\epsfbbfoundtrue
   \fi\else\ifx#1\par\else\epsffileokfalse\fi\fi}%
\def\epsfempty{}%
\def\epsfgrab #1 #2 #3 #4 #5\\{%
\global\def\epsfllx{#1}\ifx\epsfllx\epsfempty
      \epsfgrab #2 #3 #4 #5 .\\\else
   \global\def\epsflly{#2}%
   \global\def\epsfurx{#3}\global\def\epsfury{#4}\fi}%
\def\epsfsize#1#2{\epsfxsize}
\def\qed{\hspace{\fill}\hbox{${\vcenter{\vbox{                        %HOLLOW SQUARE
   \hrule height 0.4pt\hbox{\vrule width 0.4pt height 6pt
   \kern5pt\vrule width 0.4pt}\hrule height 0.4pt}}}$}}
\newcommand{\Z}{{\mathbf{Z}}} %% {{\sf Z\hspace*{-0.9ex}Z}}
\newcommand{\Q}{{\mathbf{Q}}}
\newcommand{\F}{{\mathbf{F}}}
\newcommand{\module}[2]{{{#1}\langle{#2}\rangle}}
\newcommand{\Zmodule}[1]{{\Z\langle{#1}\rangle}}
\newcommand{\Qmodule}[1]{{\Q\langle{#1}\rangle}}
\newcommand{\Omegatwelve}{{\Omega_{^{_{\leq 12}}}^{}}}
\newcommand{\field}{{k}}

\newcommand{\K}{{\cal K}}
\newcommand{\PV}{{\cal PV}}
\newcommand{\V}{{\cal V}}
\newcommand{\A}{{\cal A}}
\newcommand{\B}{{\cal B}}
\renewcommand{\P}{{\cal P}}
\newcommand{\rk}{{\rm rk\ }}
\newcommand{\image}{{\rm im\ }}
\newcommand{\iso}{{\;\;\cong\;\;}}
\newcommand{\epsfig}[1]{\noindent\epsfbox{#1}}
\newcommand{\cepsfig}[1]{
\noindent
\smallskip \par
\centerline{\epsfig{#1}}
\par
}
\newcommand{\lepsfig}[1]{
\noindent
\begin{tabular}{c}
\epsfig{#1} \\
\end{tabular}
}
\newcommand{\newln}{\hspace{\fill}\\}
\newcommand{\proof}{{\bf Proof\quad}}
\newcommand{\conjecture}{{\bf Conjecture\quad}}

\newtheorem{theorem}{Theorem}[section]
\newtheorem{definition}[theorem]{Definition}
\newtheorem{remark}[theorem]{Remark}
\newtheorem{lemma}[theorem]{Lemma}
\newtheorem{corollary}[theorem]{Corollary}
\newtheorem{proposition}[theorem]{Proposition}

\setlength{\topmargin}{-2cm}
\setlength{\textwidth}{14cm}
\setlength{\textheight}{23.5cm}
\setlength{\smallskipamount}{0.2cm}

\def\runninghead#1#2{\pagestyle{myheadings}
\markboth{{\protect\footnotesize\it{\quad #1}}\hfill}
{\hfill{\protect\footnotesize\it{#2\quad}}}}
\headsep=15pt

\begin{document}
\runninghead{Jan A.~Kneissler}{The number of primitive Vassiliev invariants up to degree twelve}

\title{The number of primitive Vassiliev invariants\\ up to degree twelve}
\author{Jan A.~Kneissler}
\date{\normalsize{\it first version: February 20, this version: June 16, 1997}}
\maketitle

\vspace{6pt}

\begin{abstract}
{
\leftskip 1cm
\rightskip 1cm
\noindent
We present algorithms giving upper and lower bounds for
the number of independent primitive rational Vassiliev invariants
of degree $m$ modulo those of degree $m-1$.
The values have been calculated for the formerly unknown degrees $m = 10, 11, 12$.
Upper and lower bounds coincide, which reveals that all Vassiliev invariants
of degree $\leq 12$ are orientation insensitive and are coming from
representations of Lie algebras so and gl.
Furthermore, a conjecture of Vogel is falsified and it is shown that the $\Lambda$-module
of connected trivalent diagrams (Chinese characters) is not free.
}
\end{abstract}

\section{Introduction}
\label{sectintro}
\subsection{Vassiliev invariants}
In the year 1990, V.~A.~Vassiliev introduced \cite{Va} a new type of knot invariants that
include the information of most of the invariants that followed the
celebrated discovery of the Jones polynomial (\cite{Jo}, \cite{H}, \cite{Ka}).

%%They can be characterized as follows.
An immersion of the circle $S^1$ in the three-sphere $S^3$ having exactly
$m$ double points and no other singularities is called $m$-singular.
Let $\K_m$ denote the set of ambient isotopy classes of $m$-singular
immersions. The elements of $K_0$ are classes of embeddings, i.e. knots in
the classical sense.

Any knot invariant $v$ with values in an abelian group,
can be extended to singular knots inductively,
by use of the desingularisation rule:
$$v(K) := v(K_+) - v(K_+)$$
where $K \in \K_m$, $K_+,K_- \in \K_{m-1}$ differ only locally like this:
\cepsfig{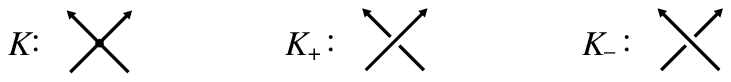}

\begin{definition}
\rm
A knot invariant with values in an abelian group is called a
{\it Vassiliev invariant} of degree $m$ iff it vanishes on $\K_{m+1}$ and the
unknot but not on $\K_m$.
A Vassiliev invariant $v$ is called {\it primitive} iff it is a monoid homomorphism,
i.e.~$v(K_1\#K_2) = v(K_1)+v(K_2)$, where $\#$ denotes the connected sum operation
for knots.
\end{definition}

\noindent
For a commutative ring $k$ let $\V_m^k$ ($\PV_m^k$) denote the $k$-module generated
by all $k$-valued (primitive) Vassiliev invariants of degree $\leq m$.

The product of two invariants $v_1,v_2$ is given pointwise by
$v_1\!\!\cdot\! v_2(K) := v_1(K)v_2(K)$ for all $K \in \K_0$.
It is not hard to show that if the degrees of $v_1, v_2$ are $m_1, m_2$, then
the degree of $v_1\!\!\cdot\! v_2$ is $m_1 + m_2$.
This establishes a graded algebra structure on $\oplus \,\V_m / \V_{m-1}$.
It is even a graded Hopf algebra (the coproduct corresponds to the connected sum operation),
which explains why we restrict ourselves to primitive Vassiliev invariants.
\begin{remark}
\label{remprim}
Every Vassiliev invariant can be expressed (uniquely up to invariants
of lower degree) as a polynomial in primitive Vassiliev invariants.
\end{remark}

\noindent
Vassiliev invariants have been defined topologically but they are
closely related to purely combinatorial objects, which we shall describe now.

\subsection{Modules of diagrams}
\begin{definition}
\rm
\begin{enumerate}
\item
A {\it free diagram} (or Chinese character) of degree $(m,u)$ is a finite abstract graph with
$2m-u$ trivalent and $u$ univalent vertices.
The trivalent vertices are rigid, i.e.~a cyclic ordering
of the three arriving edges is chosen at every trivalent vertex.
\item
A diagram together with a linear ordering of its univalent vertices
is called {\it fixed diagram}.
\item
A diagram of degree $(m,u)$ with $u > 0$ together with a cyclic ordering of its univalent vertices
is called {\it circle diagram} of degree $m$.
\end{enumerate}
\end{definition}

\noindent
In all our pictures the edges are ordered counterclockwise at each vertex.
In circle diagrams we depict the cyclic ordering of the univalent vertices
by gluing them on an oriented circle.
We will need four types of local relations (only the changed parts
of the diagrams are shown):
\begin{itemize}
\item
AS (antisymmetry of vertices): \lepsfig{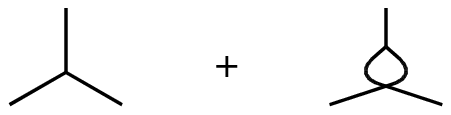}
\item
IHX relation: \lepsfig{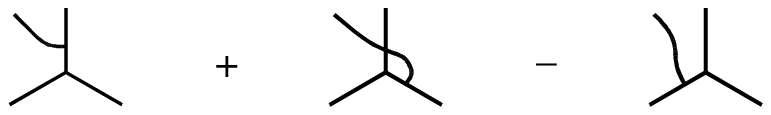}
\item
STU relation: \lepsfig{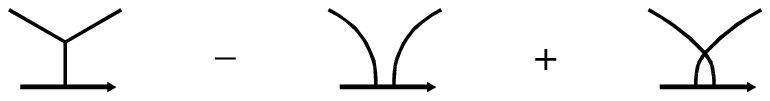}
\item
FI (framing independence): \lepsfig{}
\end{itemize}

\newln
AS and IHX are homogenous with respect to $m$ and $u$, the STU-relation only
with respect to $m$. The STU and FI relations are defined only for circle diagrams.

\begin{definition}
\rm
We have the following $\Z$-modules\footnote{ We will denote the free $k$-module with basis $S$ by $k\langle S \rangle$.}:
\begin{eqnarray*}
\A_m &:=& \Zmodule {\mbox{ circle diagrams of degree } m\; } \;/\;
\Zmodule {\mbox { STU relations }} \\
\A_m^{\rm r} &:=& \A_m / \Zmodule {\mbox { FI relations }} \\
\P_m &:=& \mbox {submodule of } \A_m \mbox{ generated by connected circle diagrams } \\
\P &:=& \bigoplus_{m=2}^\infty\limits \P_m \\
\B_{m,u} &:=& \Zmodule {\mbox{ connected free diagrams of degree } (m,u)\; } \;/\;
\Zmodule {\mbox { AS,IHX relations }} \\
F(u) &:=& \Zmodule {\mbox{ conn.~fixed diagrams with } u \mbox{ univalent vertices }} \;/\;
\Zmodule {\mbox { AS,IHX relations }} \\
\end{eqnarray*}
\end{definition}

\noindent
The most important (and highly non-trivial) facts about Vassiliev
invariants may be summarized in the following manner.
\begin{theorem}(Bar-Natan, Birman-Lin, Drinfeld, Kontsevich, Vassiliev)
\begin{eqnarray*}
&\V_m^\Q \;/\; \V_{m-1}^\Q \iso \A_m^{\rm r} \otimes \Q & \mbox { for }m \geq 1 \\
&\PV_m^\Q \;/\; \PV_{m-1}^\Q \iso \P_m \otimes \Q \iso \bigoplus_{u=1}^m\limits\; \B_{m,u}\otimes\Q &\mbox{ for } m \geq 2 \\
\end{eqnarray*}
\end{theorem}
\begin{remark}
\rm
To be more specific, there is a natural way to define a map
$\V_m^\Q / \V_{m-1}^\Q \rightarrow {\rm Hom}_\Z^{} (\A_m^{\rm r},\Q)$,
which turns out to be the desired isomorphism. $\A_m, \P_m, \B_{m,u}$ 
are finite dimensional, so we use them instead of their duals. 
\end{remark}

\noindent
It is very annoying that our knowledge of $\P$ %% = \bigoplus_{m=2}^\infty \P_m$ 
is so limited.
Dror Bar-Natan has computed $\rk \P_m$ for $m \leq 9$.
Upper and lower bounds for all degrees have been found (\cite{CD1}, \cite{CD2}, \cite{St}) but they
are unacceptably bad. And we know practically nothing about torsion in $\P$.

The main goal of this paper is to describe two algorithms that give
upper bounds for the rank of $\P_m$.
But first, we present a very good lower bound that is due
to Dror Bar-Natan and an algebra $\Lambda$, introduced by Pierre Vogel,
that acts on $\P$.

\subsection{Marked surfaces}
\label{sectmark}
If every edge of a free diagram is labeled with exactly one of the symbols "=" or "x",
it is called a {\it marked diagram}. A {\it marked surface} is a closed
compact surface with some points marked on its boundary.
At each marked point, an orientation of the boundary component is specified.
A marked surface $F$ is {\it normalized}, if either $F$ is orientable and
all markings induce the same orientation on $F$, or $F$ is non-orientable and
the orientations of the markings coincide on each component of $\partial F$.
We call a diagram {\it embedded} if it is drawn on the 2-sphere $S^2$ and
the cyclic order given at each vertex is compatible with the orientation of $S^2$.

We will thicken the five building blocks of embedded marked diagrams (univalent vertices,
trivalent vertices, edges with "=", edges with "x", crossings of edges) like this:
\cepsfig{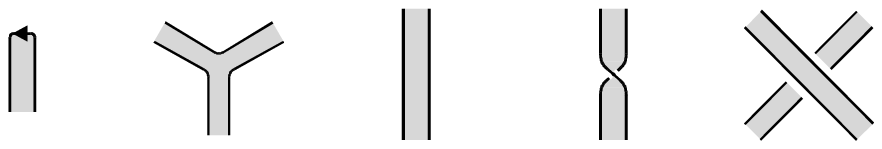}

\noindent
This assigns to every embedded marked diagram $D$ a marked surface $\hat D$.
If $D$ has $u$ univalent vertices then $\hat D$ has $u$ markings on its boundary.
If $D^\prime$ is the diagram that is obtained after all the markings of
a marked diagram $D$ are forgotten, we say that $D$ is a {\it marking of} $D^\prime$.
Let $x(D)$ denote the number of "x"-marked edges of $D$.

\noindent
Now we can define the "thickening map" $\Phi_m$ from $\bigoplus_{u=1}^{m}\B_{m,u} \rightarrow \Zmodule{\mbox{marked surfaces}}$:
\begin{definition}
\rm
If an element $b$ of $\B_{m,u}$ is represented by a embedded diagram $D_b$ then
we define $\Phi_m(b) := \sum_{{\mbox{\tiny all markings}} \atop D {\rm\; of\; } D_b}\limits (-1)^{x(D)}\hat D$.
\end{definition}

\begin{remark}
\rm
It is easy to show that $\Phi_m$ is well-defined, e.g. it does not depend on the
choice of the embedding and it respects the relations AS and IHX. This implies
that $\rk (\image\Phi_m)$ is a lower bound for $\rk \P_m$.
\end{remark}

\noindent
Let $p$ denote the projection from $\Zmodule{\mbox{marked surfaces}}$
onto $\Zmodule{\mbox{normalized marked surfaces}}$ and let $\tilde\Phi_m := p\circ\Phi_m$.
Of course we then have $\rk (\image\tilde\Phi_m) \leq \rk \P_m$ as well.
It can be shown that $\tilde\Phi_m$ contains essentially the same
information as $\Phi_m$.

\begin{remark}
\label{remlie}
\rm
On one side, to every marked surface $s$, there is naturally associated a linear
form on $\P_m$, namely the coefficient of $s$ in $\Phi_m$.
On the other side, to a finite dimensional Lie algebra, equipped with
a symmetric, Ad-invariant, non-degenerate, bilinear form and a finite
dimensional representation, there is associated a linear form on $\P_m$, too.

\noindent
Bar-Natan has shown in \cite{BN1} that all the linear forms obtained via marked surfaces
are coming from Lie algebras in the families so and gl and all of their
representations. He has also shown that the corresponding Vassiliev invariants
contain the same information as the HOMFLY and the ($2$-variable) Kauffman
polynomials and all of their cablings.
\end{remark}

\subsection{The algebra $\Lambda$}

Vogel defined an interesting submodule $\Lambda$ of the module $F(3)$ of fixed diagrams with
three univalent vertices.
The symmetric group $S_u$ acts on $F(u)$ by permutation of the univalent
vertices. There are maps $\phi_i$ ($1 \leq i \leq u$) from $F(u)$ to $F(u+1)$,
given by gluing a trivalent vertex to the $i$-th univalent vertex; the two
new univalent vertices get the numbers $i$ and $i+1$ and the numbers $> i$
are increased by one.

\begin{definition}
\rm
\label{deflambda}
$
\forall u\in F(3): \; u\in\Lambda \;\;:\Leftrightarrow\;\; \sigma(u) = \epsilon(\sigma)u
\mbox{ for all }\sigma\in S_3\;\mbox{ and }\phi_1(u) = \phi_2(u)$,
where $\epsilon$ is the signature homomorphism.
\end{definition}

\noindent
For a diagram $d$ of $\P$ and a diagram $u$ of $\Lambda$ one makes
a simple construction: Delete a trivalent vertex
of $d$ (it always has at least one) and the three univalent vertices
of $u$. The three remaining open edges of $d$ are glued to those
of $u$. The first condition in definition \ref{deflambda} together with the
AS relation cause that all $6$ ways of doing this give the result (in $\P$).
The second condition has the effect that it does not matter, at which trivalent
vertex $u$ is inserted (here it is essential that $d$ is connected).
It is easy to show that the insertion is compatible with the IHX relation,
so $\Lambda$ operates on $\P$.
$\Lambda$ is even a graded algebra because it acts on itself, and
$\P$ is a $\Lambda$-module.

\noindent
It has been shown in \cite{Vo} that $\Lambda\otimes\Q$ is commutative and that
it is contained in $\P\otimes\Q$:

\begin{proposition}
\label{vogelprop}
$\Lambda\otimes\Q \iso \bigoplus_{m=2}^\infty\B_{m,2}\otimes\Q$
\end{proposition}

\noindent
Furthermore the following elements $t,x_3,x_4,x_5,\ldots$ of $\Lambda$ 
are constructed:
\cepsfig{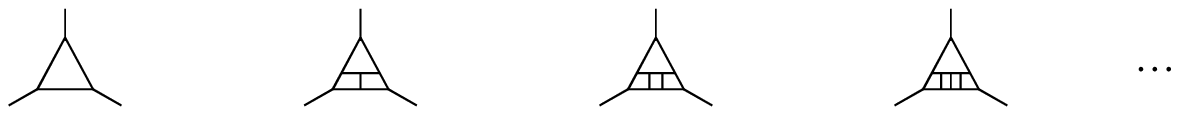}

\noindent
Vogel showed that, in degree $\leq 8$,
$\Lambda\otimes\Q$ is generated by $t,x_3,x_5,\ldots$ and isomorphic
to $\Q[t,x_3,x_5,\ldots]$. He conjectured that this is true in all degrees,
and gave a polynomial in degree $10$ for which he could not tell,
whether it is trivial or not.

\section{Results}
\label{sectres}

We have implemented both algorithms that will be given in section \ref{sectalg} and made
a program that effectively computes the thickening map $\tilde\Phi_m$ described in section
\ref{sectmark}.
The output of these three programs for degree $m$ will be denoted $O_A(m),
O_B(m),O_C(m)$,
respectively. By corollary \ref{corineq} we have that $O_A(m),O_B(m) \geq \rk \P_m \geq O_C(m)$.
First we confirmed for $3 \leq m \leq 9$ that $O_A(m) = O_B(m) = O_C(m) =$ values given in \cite{BN2}.
Then we found $O_A (10) = O_B(10) = O_C(10) = 27$, $O_A(11) = O_B(11) = O_C(11) = 39$.
We were astonished that our algorithms are good enough to give
the exact values. The complete surprise came with $m = 12$;
we had conjectured $O_C(12) \neq \rk \P_{12}$, because of Vogel's conjecture.
After some hundred hours of CPU-time we got the result $O_B(12) = O_C(12) = 55$
(the computation of $O_A(12)$ is too exhaustive and could not be performed).

\subsection{Vassiliev invariants}
\begin{theorem}[computational result]
The sequence of the dimensions of the spaces of rational-valued primitive Vassiliev invariants of degree $m$
modulo those of degree $m-1$ starts
$$0,\;1,\;1,\;2,\;3,\;5,\;8,\;12,\;18,\;27,\;39,\;55,\;\ldots$$
\end{theorem}

\begin{remark}
{\rm
For some time only the first seven numbers of this sequence were known and there
was some excitement, because it appeared to be the famous Fibonacci
sequence. It is a somehow mysterious coincidence that it
is again a Fibonacci number in degree twelve.
}
\end{remark}

\noindent
In view of remark \ref{remlie} we have the following consequence of $O_C(m) = \rk \P_m$
for $m \leq 12$.

\begin{corollary}
All rational Vassiliev invariants of degree $\leq 12$ are coming from representations
of the classical Lie algebras so and gl.
\end{corollary}

\begin{remark}
{\rm
Vogel has shown that a similar statement is false for sufficiently large
degrees. Let $m_c$ denote the minimal degree for which not
all Vassiliev invariants are coming from semi-simple Lie algebras.
Results of Jens Lieberum (\cite{Li}) and our calculations together imply that
$13 \leq m_c \leq 17$.}
\end{remark}

\begin{corollary}
Vassiliev invariants up to degree twelve can not distinguish knots from
their inverses.
\end{corollary}
\proof
It has been shown that Vassiliev invariants coming from semi-simple
Lie algebras are orientation insensitive. Another way is to
verify that $\rk \B_{m,u} = 0$ for $m \leq 12$ and $u$ odd (see table
in section \ref{secttab}).
\qed

\newln
A statement about ${\cal A}$ similar to remark \ref{remprim} makes it is
easy to compute the ranks of $A_m^{\rm r}$ and $A_m$,
corresponding to the numbers of Vassiliev invariants of knots
and Vassiliev invariants of framed knots, if
$\rk \P_m$ is known.

\newln
\centerline{\begin{tabular}{||c|c|c|c|c|c|c|c|c|c|c|c|c|c||}
\hline
$m$ &
	$0$ &   $1$ &   $2$ &   $3$ &   $4$ &   $5$ &   $6$ &   $7$ &
	$8$ &   $9$ &   $10$ &  $11$ &   $12$  \\
\hline
$\rk\P_m$ &
	$0$ &   $1$ &   $1$ &   $1$ &   $2$ &   $3$ &   $5$ &   $8$ &
	$12$ &  $18$ &  $27$ &  $39$ &   $55$ \\
$\rk\A_m$ &
	$1$ &   $1$ &   $2$ &   $3$ &   $6$ &   $10$ &   $19$ &   $33$ &
	$60$ &  $104$ & $184$ &  $316$ &   $548$ \\
$\rk\A_m^{\rm r}$ &
	$1$ &   $0$ &   $1$ &   $1$ &   $3$ &   $4$ &   $9$ &   $14$ &
	$27$ &  $44$ &  $80$ &  $132$ &   $232$   \\
\hline
\end{tabular}}

\subsection{Results about $\P$ and $\Lambda$}
\label{secttab}

Our programs work over the field $\F_2$, so due to corollary \ref{cortorsion}
we have a little statement about torsion in $\P$.

\begin{corollary}
There is no $2$-torsion in $\P$ in degree $\leq 12$.
\end{corollary}

\noindent
By counting the dimensions of the image of the thickening map $\Phi$ for
each $\B_{m,u}$ separately, we get the following table for $\rk \B_{m,u}$.

\newln
\centerline{\begin{tabular}{||c|c|c|c|c|c|c|c|c||}
\hline
$\rk \B_{m,u}$&       $u=2$ & $u=4$ & $u=6$ & $u=8$ & $u=10$ & $u=12$ & total \\
\hline
$m=1$ & $1$ &   &       &       &       &       &       $1$ \\
\hline
$m=2$ & $1$ &   &       &       &       &       &       $1$ \\
\hline
$m=3$ & $1$ &   &       &       &       &       &       $1$ \\
\hline
$m=4$ & $1$ &   $1$ &   &       &       &       &       $2$ \\
\hline
$m=5$ & $2$ &   $1$ &   &       &       &       &       $3$ \\
\hline
$m=6$ & $2$ &   $2$ &   $1$ &   &       &       &       $5$ \\
\hline
$m=7$ & $3$ &   $3$ &   $2$ &   &       &       &       $8$ \\
\hline
$m=8$ & $4$ &   $4$ &   $3$ &   $1$ &   &       &       $12$ \\
\hline
$m=9$ & $5$ &   $6$ &   $5$ &   $2$ &   &       &       $18$ \\
\hline
$m=10$ & $6$ &   $8$ &   $8$ &   4 &  $1$ &     &       $27$ \\
\hline
$m=11$ & $8$ &   $10$ &  $11$ &  $8$ & $2$ &     &      $39$ \\
\hline
$m=12$ & $9$ &   $13$ &  $15$ &  $12$ & $5$ & $1$ &     $55$ \\
\hline
\end{tabular}}

\begin{corollary}
The algebra morphism from $\Q[T,X_3,X_5,\ldots]$ to
$\Lambda\otimes\Q$ given by $T \rightarrow t$, $X_i \rightarrow x_i$ is
not an isomorphism. In degree $< 11$ it is surjective and has a one dimensional kernel
(living in degree $10$).
\end{corollary}
\proof
By calculating characters, Vogel has already shown that this algebra
morphism is injective in degree $\leq 9$ and its kernel in degree $10$ has at
most dimension one. Because of proposition \ref{vogelprop}, the statements can be verified
by counting the number of monomials in each degree $m$ and comparing it
to $\rk \B_{m+2,2}$.
In degree $10$, for example, there are ten monomials ($t^{10},t^7x_3^{},t^5x_5^{},t^3x_7^{},
tx_9^{},t^4x_3^2,t^2x_3^{}x_5^{},x_3^{}x_7^{},x_5^2$ and $tx_3^3$), but we
have $\rk \B_{12,2} = 9$.
\qed

\newln
So one half of Vogel's conjecture is false, but the calculations show that up
to degree $10$ the other one holds:

\newln
\conjecture
$\Lambda\otimes\Q$ is generated (as algebra over $\Q$) by the elements $t,x_3,x_5,\ldots$

\begin{corollary}
$\P\otimes\Q$ is not a free $\Lambda\otimes\Q$-module.
\end{corollary}
\proof
Let us assume that $\B_{m,4}\otimes\Q$ is a free $\Lambda\otimes\Q$-module
with rank $\alpha_m \geq 0$ ($m \geq 4$). Let $\lambda_m$ denote the dimension
of $\Lambda\otimes\Q$ in degree $m$. Then we have the following
formula for the rank of $\B_{m,4}$:
$$\rk \B_{m,4} = \dim_\Q \B_{m,4}\otimes\Q = \sum_{i=4}^{m} \lambda_{m-i}\alpha_{i}$$
We have $\lambda_0,\ldots,\lambda_7 = 1,1,1,2,2,3,4,5$,
which together with the values $\rk B_{4,4},\ldots,\rk \B_{11,4}$
implies $\alpha_4 = \alpha_6 = \alpha_8 = \alpha_{10} = 1$,
$\alpha_5 = \alpha_7 = \alpha_9 = 0$, $\alpha_{11} = -1$.

This contradiction shows that $-$ at least in the $u = 4$ column
of $\P\otimes\Q$ $-$ nontrivial relations hold.\qed

\newln
We have found another relation which is located in the $u=6$ column in degree $12$.
Unlike here, its existence can not be shown by simply counting dimensions.

\subsection{The structure of $\P\otimes\Q$ as far as we know it}

Using the thickening map, we have built a minimal set of diagrams $\Omegatwelve$
that generate $\P\otimes\Q$ as $\Lambda\otimes\Q$-module in degrees up to $12$.
We were trying to make the elements of $\Omegatwelve$ as simple as possible and
finally, this lead us to a very special type of diagrams:
\begin{definition}
Let $\omega_{i_1i_2\ldots i_k}$ denote the element of $\B_{i_1+\ldots+i_k+k-1,\;i_1+\ldots+i_k}$
that is represented by a "caterpillar" diagram consisting
of $k$ "body segments" with $i_1,\ldots,i_k$ "legs", respectively.
\end{definition}
Here are some examples of caterpillar diagrams for $\omega_4,\;\omega_{302},\;\omega_{13},\;\omega_{02131}$.
\cepsfig{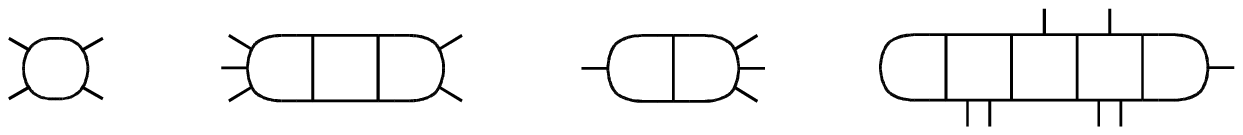}

\begin{remark}
\label{remcat}
\rm
It is a nice exercise to use the AS and IHX relations to prove that $\omega_{i_1i_2\ldots i_k}$
is well defined (i.e. for inner segments it makes no difference on which
side of the body the legs are drawn).
An easy consequence is $\omega_{i_1\ldots i_k} = \omega_{i_k\ldots i_1}$.
The diagrams $\omega_i$ are also called "wheels with $i$ spokes".
\end{remark}

\noindent
Let $\Omegatwelve$ denote the set consisting of the following $31$ elements.
\newcommand{\w}[1]{\omega_{#1}^{}}
\newcommand{\oneoftwo}[1]{\vbox{\hbox{$#1$}\vspace{3pt}}}
\newcommand{\twooftwo}[2]{\vbox{\vspace{3pt}\hbox{$#1$}\vspace{3pt}\hbox{$#2$}}}
\newcommand{\threeofthree}[3]{\vbox{\vspace{3pt}\hbox{$#1$}\vspace{3pt}\hbox{$#2$}\vspace{3pt}\hbox{$#3$}}}
\newcommand{\oneofthree}[1]{\vbox{\hbox{$#1$}\vspace{7pt}}}
\newcommand{\twoofthree}[2]{\vbox{\hbox{$#1$}\vspace{3pt}\hbox{$#2$}\vspace{3pt}}}
\noindent
\begin{center}
\begin{tabular}{||c|c|c|c|c|c|c||}
\hline
$\Omegatwelve_{_{}}$&$ u=2 $&$ u=4 $&$ u=6 $&$ u=8 $&$ u=10 $&$ u=12 $\\
\hline
$m=2 $&$ \w{2} $&$  $&$  $&$  $&$  $&$  $\\
\hline
$ m=3 $&$  $&$  $&$  $&$  $&$  $&$  $\\
\hline
$ m=4 $&$  $&$ \w{4} $&$  $&$  $&$  $&$  $\\
\hline
$ m=5 $&$  $&$  $&$  $&$  $&$  $&$  $\\
\hline
$ m=6 $&$  $&$ \w{202} $&$ \w{6} $&$  $&$  $&$  $\\
\hline
$ m=7 $&$  $&$  $&$ \w{24} $&$  $&$  $&$  $\\
\hline
$ m=8 $&$  $&$ \w{20002} $&$ \w{204} $&$ \w{8} $&$  $&$  $\\
\hline
$ m=9 $&$  $&$  $&$ \w{2004} $&$ \w{26} $&$  $&$  $\\
\hline
$\oneoftwo{\,\,\,m=10} $&$  $&$ \oneoftwo{\w{2000002}} $&$ \twooftwo{\w{20004}}{\w{20202}} $&$
\twooftwo{\w{206}}{\w{404}} $&$ \oneoftwo{\w{1\!0}}  $&$  $\\
\hline
$\oneofthree{\,\,\,m=11} $&$  $&$  $&$ \oneofthree{\w{200004}} $&$ \threeofthree{\w{2006}}{\w{4004}}{\w{2222}} $&$ \oneofthree{\w{28}}  $&$  $\\
\hline
$\oneofthree{\,\,\,m=12} $&$  $&$ \oneofthree{\w{200000002}} $&$ \twoofthree{\w{2000004}}{\w{2000202}} $&$ \threeofthree{\w{20006}}{\w{40004}}{\w{20204}} $&$ \threeofthree{\w{208}}{\w{406}}{\w{262}} $&$ \oneofthree{\w{1\!2}} $\\
\hline
\end{tabular}
\end{center}

\begin{remark}
\rm
At several places in the upper table the choice of a minimal generating set is not
unique. We tried to make $\Omegatwelve$ look as uniform as possible. At first
place we were able to renounce on $\omega$'s with odd indices. After this
only few choices still had to be done. For example, 
we preferred $\omega_{26}$ over $\omega_{44}$ (because of the other entries
of the form $\omega_{2,u-2}$) and $\omega_{2222}$ over $\omega_{2204}$ (because
of its symmetry).
\end{remark}

\noindent
Let $\P_\omega$ denote sub-$\Lambda$-module of $\P$ that is generated by caterpillar diagrams.
A glance through the table $\Omegatwelve$ immediately opens the following two questions:

\begin{enumerate}
\item
Is $\P_\omega$ already generated by the caterpillar diagrams with even indices?
\item
Is $\P_\omega = \P$?
\end{enumerate}
The AS relation and remark \ref{remcat} implies that caterpillar diagrams with an odd number of univalent vertices are always trivial.
So if question 2.~could be answered positively, it would imply that all
Vassiliev invariants are orientation insensitive.

It is tempting to make conjectures about how this table continues (especially
for the $u = 4$ column), but let us just summarize what we know for sure.

\begin{remark}
\label{remom}
\rm
If $\Omega = \bigcup \Omega_{m,u}$ is a minimal set of free diagrams that
generate $\P\otimes\Q$ as $\Lambda\otimes\Q$-module, then
\begin{itemize}
\item
The $u = 2$ column is essentially empty: $\Omega_{m,2}
= \emptyset$ for $m > 2$.
\item
The $u = m = 2i$ diagonal consists only of wheels: $\Omega_{2i,2i} = \{\;\omega_{2i}\;\}$.
\item
On the first subdiagonal ($m-1 = u = 2i$) we have exactly $\lfloor {i\over 3}\rfloor$ elements.
A natural choice is $\Omega_{2i+1,2i} := \{ \;\omega_{ab} \;\vert\; a > 0;\; a$ even$;\; 2a\leq b;\; a+b = 2i\;\}$.
\item
On the second subdiagonal ($m-2 = u = 2i$) there are exactly $\#\Omega_{2i+2,2i} = \lfloor{(i+1)^2\over 12}+{1 \over 2}\rfloor$ elements.
\item
For odd $u$ we know $\Omega_{m,u} = \emptyset$ if $m-u \leq 5$ or $u = 1$ or $m \leq 12$.
\end{itemize}
\end{remark}
\proof
The first statement is due to proposition \ref{vogelprop}. The second is obvious and the third
and forth follow from results of Oliver Dasbach in \cite{Da2}.
He showed that $\Omega_{2i+1,2i}\cup\{\;t\omega_{2i}\;\}$ is a basis
for $\B_{2i+1,2i}$ and $\dim\B_{u+2,u}=\lfloor { {u^2+12u} \over {48} }\rfloor+1$
for $u$ even.
It can be verified that $\{ \;t^2\omega_{2i}\;\}\cup \; t\,\Omega_{2i+1,2i}$
are independent in $\B_{2i+2,2i}$, so $\#\Omega_{2i+2,2i} = \lfloor { {4i^2+24i} \over {48} } \rfloor + 1
- \lfloor {i\over 3}\rfloor - 1 = \lfloor{(i+1)^2\over 12}+{1 \over 2}\rfloor$. 
Dasbach has also shown (\cite{Da1}) that $\B_{m,u} = 0$ for $m-u \leq 5$ and $u$ odd.
The same statement for $u=1$ is easy to prove.
\qed

\begin{remark}
\rm
Some time after having found $\Omegatwelve$, we discovered that caterpillar
diagrams have already been used by Chmutov and Duzhin in \cite{CD2}, who call them
"baguette diagrams". The main theorem of \cite{CD2} states that the elements $\omega_{n_1\ldots n_k}$
with all $n_i$ even, $\sum_{i=1}^{j-1} n_i < n_j$ for $j < k$ and
$\sum_{i=1}^{k-1} n_i < {1 \over 2}n_k$ are linearly independent.
This result is quite striking, but it is useless in our context, because
the first interesting\footnote{ The third item of remark \ref{remom} and
the simple relation $\omega_{0,n_2\ldots n_k} = 2t\omega_{n_2\ldots n_k}$ should make clear
that only the diagrams $\omega_{n_1\ldots n_k}$ with $n_1 > 0$ and $k > 2$
are of further interest.}
diagram in this set is $\omega_{2,4,14}$ and lies in
degree $22$.
\end{remark}

\renewcommand{\w}[1]{w_{#1}^{}}
\newcommand{\x}[2]{X_{#1}^{#2}}
\renewcommand{\t}[1]{T_{}^{#1}}
\newcommand{\p}{+}
\newcommand{\m}{-}
\noindent
Let $W$ denote the image of $\Omegatwelve$ under the map $\omega_* \rightarrow w_*$
and make the following abbreviations:
{
\scriptsize
\begin{eqnarray*}
P &\!:=\!&\! 32\t{10}\! \m 152\t{7}\!\x{3}{} \p 252\t{5}\!\x{5}{}
         \m 101\t{4}\!\x{3}{2} \m 36\t{3}\!\x{7}{}  \m 9\t{2}\!\x{3}{}\x{5}{}
         \p 14\t{}\x{3}{3} \p 9\x{3}{}\x{7}{} \m 9\x{5}{2} \\
Q & := & \left(32\t{7}\! \p 35\t{4}\!\x{3}{} \m 9\t{2}\!\x{5}{} \m 4\t{}\x{3}{2}\right)\!\otimes\!\w{4}
       + \left(\m 76\t{5}\! \p 10\t{2}\!\x{3}{} \p 3\x{5}{}\right)\!\otimes\!\w{202}
       + \left(12\t{3}\! \m 3\x{3}{}\right)\otimes\w{20002} \\
R & := & \left(\m 16\t{6}\! \p 21\t{3}\!\x{3}{} \m 3\t{}\x{5}{} \m 2\x{3}{2}\right)\!\otimes\!\w{6}
       + \left(32\t{5}\! \m 17\t{2}\!\x{3}{} \p 3\x{5}{}\right)\!\otimes\!\w{24}
       + \left(\m 36\t{4}\! \p 9\t{}\x{3}{} \right)\!\otimes\!\w{204}\\
&&     + \left(12\t{3}\! \m 3\x{3}{}\right)\!\otimes\!\w{2004} \\
\end{eqnarray*}
}
\begin{theorem}
\label{theop}
There is a module morphism
$$\mu \;\;:\;\; Q[T,X_3,X_5,\ldots]\otimes\Qmodule{W} \;\;/\;\;
\Qmodule{P\otimes w_2, Q, TQ, R\;} \;\;\;\rightarrow\;\;\; \P\otimes\Q.$$
It is an isomorphism in degree $\leq 12$.
\end{theorem}
\proof
$\mu$ is given by $\mu(T^aX_3^bX_5^c\ldots\otimes w_*) = t^ax_3^bx_5^c\ldots \omega_*$.
We have seen that Bar-Natan's thickening map $\Phi$ is injective on $\P_m\otimes\Q$
for $m \leq 12$.
So we have to calculate $\Phi$ for all diagrams
of the form $t^ax_3^bx_5^c\ldots \omega_*$ ($\omega_* \in \Omegatwelve$) with
degree $\leq 12$ (there are exactly $175$ of them).
The program and the results are available via the internet
(see section \ref{hintsec}).
It is then easy to verify that the span of the resulting vectors
is $171$ dimensional and that
$P\otimes w_2$, $Q$, $TQ$ and $R$ span the kernel of $\Phi\circ\mu$.
So $\mu$ is well defined, injective, and because of
$\sum_{m=2}^{12}\nolimits \rk \P_m = 171$ it is also surjective. \qed

\section{The principle behind the algorithms}
\label{sectprin}

We will now describe a prototype of an algorithm that yields an upper
bound for the rank (ubr) of a finitely generated abelian group $A$.

\begin{definition}
{\rm
A quintuple $(k,S,\varphi,\delta,\rho)$ where $k$ is a field,
$S$ is a finite set, $\varphi$ is a mapping $\varphi: S \rightarrow A$
and $\delta,\rho$ are endomorphisms of $\module kS$ shall be called
{\it ubr-algorithm} for the finitely generated abelian group $A$, iff the following conditions are satisfied:
{
\setlength{\parskip}{-3pt}
\begin{enumerate}
\setlength{\itemsep}{-2pt}
\item
$\varphi(S)$ is a set of generators of $A$,
\item
there exists a integer $j$ such that $\delta^{j+1} = \delta^{j}$,
\item
$\hat\varphi \circ \delta = \hat\varphi$,
\item
$\hat\varphi \circ \rho = 0$.
\end{enumerate}
}
\noindent
Here $\hat\varphi$ denotes the vectorspace homomorphism
$\hat\varphi: \module kS \rightarrow A\otimes_\Z^{}\!k$ that is induced by $\varphi$.
}
\end{definition}

\begin{definition}
\label{defubralg}
{\rm
For a given ubr-algorithm let $\Delta := \delta^j$ and ${\cal I} := {\rm im}\, \Delta$.
Then $\bar\rho := \Delta\circ\rho\vert_{{\cal I}}$ is an ${\cal I}$-endomorphism.
The {\it output} of the ubr-algorithm defined as the natural number
${\rm output}(k,S,\varphi,\delta,\rho) := \dim_k \ker \bar\rho$.
}\end{definition}

\begin{lemma}
\label{helplemma}
Conditions 2 and 3 imply
$\Delta(\ker \hat\varphi) = \ker \hat\varphi \cap {\cal I}$.
\end{lemma}
{\noindent\bf Proof:}
"$\supset$": 2.$\;\Rightarrow\; \Delta$ is a projection onto ${\cal I}$.
\\
$"\subset"$: 3.$\;\Rightarrow\; \hat\varphi \circ\Delta = \hat\varphi
\;\Rightarrow\; \Delta(\ker \hat\varphi) \subset \ker \hat\varphi$.
\qed

\begin{proposition}
If $(k,S,\varphi,\delta,\rho)$ is an ubr-algorithm for the finitely generated
abelian group $A$ then
${\rm output}(k,S,\varphi,\delta,\rho) \geq \rk A$.
\end{proposition}

{\noindent\bf Proof:}
The first condition of definition \ref{defubralg} implies that $\hat\varphi$ is an epimorphism.
Due to 4.~we have $\rho({\cal I}) \subset \ker \hat\varphi$.
Together with lemma \ref{helplemma} and $\hat\varphi \circ\Delta = \hat\varphi$
we get:
\begin{eqnarray*}
\dim_k (A\otimes_\Z^{}\! k) &
 = & \dim_k \hat\varphi(\module kS) \;\;=\;\; \dim_k \hat\varphi\circ\Delta (\module kS) \;\;=\;\; \dim_k \hat\varphi({\cal I})\\
&= &\dim_k {\cal I} - \dim_k \ker (\hat\varphi\vert_{{\cal I}})
 \;\;=\;\; \dim_k {\cal I} - \dim_k \Delta (\ker \hat\varphi)\\
&\leq& \dim_k {\cal I} - \dim_k \Delta \big(\rho ({\cal I})\big)
 \;\;=\;\ \dim_k\ker\bar\rho
\end{eqnarray*}
Let $T$ denote the maximum torsion subgroup of $A$. $A$ is isomorphic
to $T\times \Z^{\rk A}$, and thus $\dim_k (A\otimes_\Z^{}\! k) \geq \rk A$.
\qed

\begin{corollary}
\label{cortorsion}
If $(F_p,S,\varphi,\delta,\rho)$ is an ubr-algorithm for $A$ and
${\rm output}(F_p,S,\varphi,\delta,\rho) = \rk\, A$ then $A$ has no
elements of order $p$.
\end{corollary}

\begin{remark}
\rm
To calculate ${\rm output}(k,S,\varphi,\delta,\rho)$ one has to find a basis of
${\cal I}$ and to evaluate the nullity of the corresponding matrix for $\bar\rho$.
So the costs (in time and space) of an ubr-algorithm mainly depend
on the dimension of $\cal I$, whereas the quality of the result (the sharpness
of the upper bound) depends on the choice of $\rho$.
\end{remark}

\section{Two algorithms for $\P_m$}
\label{sectalg}
\subsection{Circle diagrams without loops}
\label{sectcirc}
Let $S_n$ denote the set of permutations of $n$ elements.
We will use the standard linear ordering for the permutations:
$$\pi < \phi \;\; :\Leftrightarrow \;\; \exists \, i \in \{1,\ldots,n\}: \pi(i) < \phi(i)
\mbox { and } \pi(j) = \phi(j) \mbox{ for all } 1 \leq j < i.$$
Denote by $\tau_i$ the elementary transpositions $(i \;\; i\!+\!1)$ and let the product of
permutations be defined by $(\pi\phi)(i) := \phi(\pi(i))$.

\begin{lemma}
\label{smallerlemma}
For any $\pi \in S_n$ the following two statements hold:
If there exists an integer $i$ $(1 \leq i < n)$ with $\pi(i) > \pi(i+1)$  then
$\tau_i\pi < \pi$.
If there exist $i,j$ $(1 \leq i < j \leq n)$ with $\pi(i) = \pi(j)+1$ then $\pi\tau_{\pi(j)} < \pi$.
\end{lemma}
{\bf Proof:} Simply identify the $i$ in the upper definition of $<$ with this $i$ here.
\qed

\newln
\noindent
The {\it picture} of a permutation $\pi \in S_n$ is given by two vertical lines with
$n$ distinct points marked on each, together with $n$ lines connecting the $i$-th and the
$\pi(i)$-th point (counted upwards). For example the picture of $(1\;2\;4)(3)$ is
\cepsfig{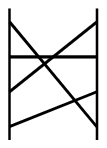}

\noindent
\begin{definition}
{\rm
For any $\pi \in S_n$ we will denote by $D^A_\pi$ the element of ${\cal P}_{n+1}$
that is obtained by replacing the box in the following figure
by the picture of $\pi$.
\cepsfig{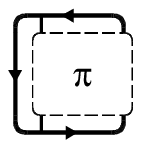}

\noindent 
The map $\varphi_A : S_n \rightarrow  \P_{n+1}$ is given by 
$\pi \rightarrow D^A_\pi$.
}\end{definition}

\noindent
We will now introduce three types of moves, which replace a permutation by
a linear combination of permutations, by showing the parts of their pictures
that are concerned. Omitted parts are indicated by dots and are assumed to
be identical in a row.

\newln
\noindent
Type I
\smallskip
\newln
\epsfig{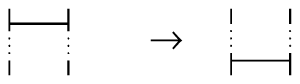}

\newln
\noindent
Type II
\smallskip
\newln
\epsfig{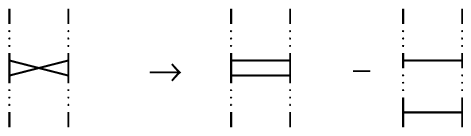}

\newln
\noindent
Type III
\smallskip
\newln
\epsfxsize = \textwidth \epsfig{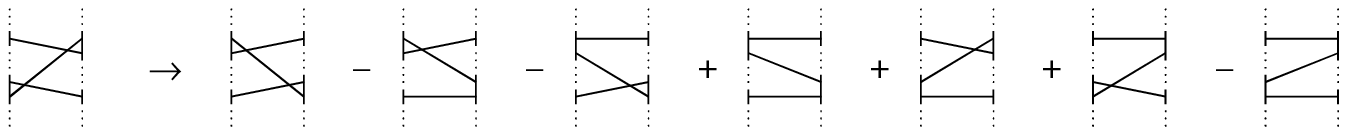}
\newln
\epsfxsize = \textwidth \epsfig{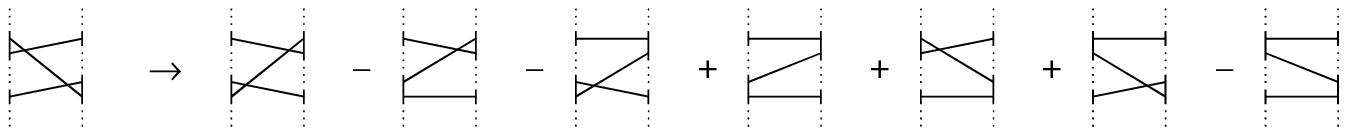}
\newln
\epsfxsize = \textwidth \epsfig{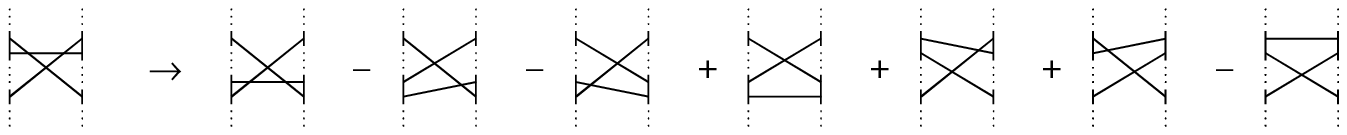}
\newln
\epsfxsize = \textwidth \epsfig{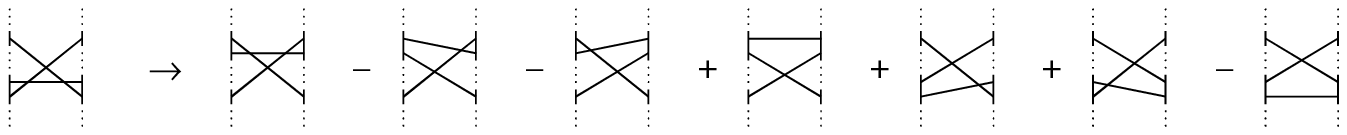}

\newln
The move I should be interpreted in the following way: when $\pi(n) = n$ then replace $\pi$
by $\tau_1\cdots\tau_{n-1}\,\pi\,\tau_{n-1}\cdots\tau_1$; the move II reads:
when there is some $i \in \{1,\ldots,n-1\}$ with $\pi(i) = \pi(i+1)+1$ then
replace $\pi$ by $\tau_i\,\pi - \tau_1\cdots\tau_{i-1}\,\pi\,\tau_{\pi(i)-1}\cdots\tau_1$; etc.

We call a move {\it reducing}, if all the permutations on the right side are
smaller than the permutation $\pi$ on the left, and we then say that $\pi$ allows
a reduction. Permutations that do not allow any reduction are called irreducible.

If a diagram allows more than one reducing move, we have to make a choice.
To make our calculations reproducible, we describe how our implementation works.
Moves of type II are indexed by the left height of the two crossing
lines that are concerned. Moves of type III are indexed by the height
of the left endpoint of the line that is "switched" during the move.
Now all possible moves of a permutation can be ordered by I $<$ II $<$ III
and II$_i$ $<$ II$_j$ and III$_i$ $<$ III$_j$ for $i < j$.

\begin{definition}
{\rm
\label{defdeltaa}
When $\pi \in S_n$ allows reductions of type I, II or III then set $\delta_A (\pi) := $ result of the
smallest reducing move.
When $\pi$ is irreducible set $\delta_A (\pi) := \pi$.
}\end{definition}

\begin{remark}
\rm
Lemma \ref{smallerlemma} implies that moves I and II are always reducing,
and that for moves of type III it suffices to check the first term on the left right.
\end{remark}

\noindent
In the following, we will work with a special element $\Theta_k$ of $\Z[S_n]$:
$$\Theta_k := \prod\limits_{i=1}^{k-1}(1-\prod\limits_{j=1}^{k-i}\tau_{k-j}) =
(1-\tau_{k-1}\tau_{k-2}\cdots\tau_1)(1-\tau_{k-1}\cdots\tau_2)\cdots(1-\tau_{k-1})$$
\begin{definition}
{\rm
\label{defrhoa}
For $\pi\in S_n$ set $k := n+1-\pi^{-1}(1)$ and let $\pi^\prime$ be
given by
$$\pi^\prime(i) = \left\{
\begin{array}{lcl}
\pi(n+1-k+i)-1 &\mbox{ for } & 1 \leq i < k \\
n &\mbox{ for } & i = k \\
\pi(n+1-i)-1 &\mbox{ for } & k < i \leq n\\
\end{array}
\right\} \mbox { (thus }\pi^\prime \in S_n\mbox{).}
$$
Now let $\rho_A$ be the $\Zmodule {S_n}$-endomorphism given via
$\rho_A(\pi) := \pi-(-1)^{n-k}\Theta_k \pi^\prime \in {\cal S}_n$.
}
\end{definition}

\noindent
This completes the description of the first algorithm 
$A_m := (\field,S_{m-1},\varphi_A,\delta_A,\rho_A)$.

\subsection{Circle diagrams with one loop}
\label{sectone}

To compare permutations of different symmetric groups, we
extend the definition of "$<$" in the following way: $i < j$, $\pi_1 \in S_i$,
$\pi_2 \in S_j \Rightarrow \pi_1 < \pi_2$.

\begin{definition}
{\rm
For $\pi \in S_n$ let $D^B_\pi$ be the element of ${\cal P}_{n}$ that is obtained by replacing
the box in the following figure by the picture of $\pi$.
\cepsfig{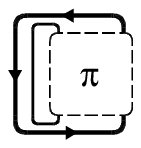}

\noindent
Let $\varphi_{B,m}$ denote mapping $\bigcup_{n=3}^{m} S_n \rightarrow \P_m$ 
that is given by $S_n \ni \pi \rightarrow t^{m-n}D^B_\pi$.
}
\end{definition}

\noindent
To define $\delta_B$ we have to modify the first two moves.
The moves of type III given in section \ref{sectcirc} 
and the ordering of moves can be adopted unchanged.

\newln
Let $\mu_n,\nu_n \in S_n$ be permutations given by $\mu_n(k)=n+1-k$ and
$\nu_n = (1\;\;2\;\ldots\;n)$.
$G := \Z/n\Z \times \Z/2\Z \times \Z/n\Z$ acts on $S_n$ by $(a,b,c)\pi := \nu_n^a\mu_n^b\pi\nu_n^c$
for all $\pi \in S_n$.
For $\pi\in S_n$ let $(\alpha,\beta,\gamma)$ denote an element of $G$ such that
$(\alpha,\beta,\gamma)\pi$ is minimal in the orbit $G\pi$ of $\pi$.
The move I$^\prime$ is to replace $\pi$ by $(-1)^{n\beta}(\alpha,\beta,\gamma)\pi$.

\newln
The move II$^\prime$ is given graphically:
\cepsfig{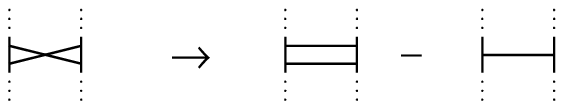}

\noindent
If the move II$^\prime$ is applied to an element of $S_n$, the second term
in the result lies in $S_{n-1}$. We will use this move only for $n \geq 4$, so
we do not have to deal with the symmetric groups $S_1$ or $S_2$.

\begin{definition}
{\rm
\label{defdeltab}
When $\pi \in S_n$ allows reductions of type I$^\prime$, II$^\prime$ or III then set $\delta_B (\pi) := $ result of the
smallest reducing move.
When $\pi$ is irreducible let $\delta_B (\pi) := \pi$.
}\end{definition}

\noindent
For any $\pi \in S_{n-1}$ and $0 \leq k \leq n-3$ we have the following
element of $\Z[S_n]:$
$$\Upsilon_{n,k}^{} (\pi) := \chi_{1,n-k-2}^{} \; \Theta_{n-k-2}^{} \; \chi_{n-k-2,k+1}^{} \; (1 - \tau_n^{})\; \pi^\#$$
Here $\Theta_{k}^{}$ is the same as in section \ref{sectcirc} and
$\chi_{r,s}^{},\pi^\# \in S_n$ are given by
$$\chi_{r,s}^{}(i) := \left\{
\begin{array}{ll}
i\!+\!s & \mbox{if } i \leq r \\
i\!-\!r & \mbox{if } r < i \leq r\!+\!s \\
i   & \mbox{if } i > r\!+\!s \\
\end{array}
\right\} \mbox{ and }
\pi^\#(i) := \left\{
\begin{array}{ll}
\pi(i) & \mbox{if } \pi(i) < \pi(n\!-\!1) \\
\pi(i)\!+\!1 & \mbox{if } \pi(i) > \pi(n\!-\!1) \\
\pi(n\!-\!1)\!+\!1 & \mbox{if } i = n \\
\end{array}
\right\}.$$

\begin{definition}
\label{defrhob}
\rm
For $\pi \in S_n$ with $\pi(1) = 1, \pi \neq {\rm id}_n$ set
$p := \max \{ \;j \;\vert \; \pi(i) = i \mbox { for all } i \leq j\;\}$
and $q := \pi^{-1}(p+1)$. If $\pi(n) \neq 2$ and $\pi(n) \neq n$ let $\pi^\prime \in S_{n-1}$
be the permutation given by
$$\pi^\prime (i) :=\left\{
\begin{array}{ll}
i & \mbox{if } i = 1 \mbox{ or } n-p < i \leq n-1\\
\pi(q+1-i) & \mbox{if } 2 \leq i \leq q-p\\
\pi(i+p) & \mbox{if } q-p < i \leq n-p \\
\end{array}
\right\}.$$
Finally for any $\pi \in S_n$ with $3 \leq n \leq m$ set
$$\rho_{B,m}^{} (\pi) := \left\{
\begin{array}{ll}
(1)(2)(3\;4\;5)+{\rm id}_5-2{\rm id}_4+{\rm id}_3 & \mbox{ if } \pi = {\rm id}_m \mbox{ and } m \geq 5 \\
\pi - \pi\tau_p +(-1)^{q-p}\Upsilon_{n,q-p-1}(\pi^\prime)& \mbox{ if }n = m \mbox{ and } \pi(1)=1 \mbox { and } \pi(n) \neq 2,n\\
\pi - \Upsilon_{n,0}(\pi)& \mbox{ if } 3 \leq n < m \mbox{ and }\pi(1)=1\\
0 & \mbox{ otherwise}
\end{array}
\right.
$$
\end{definition}

\noindent
Now we have our second candidate for an ubr-algorithm $B_m := (\field,\bigcup_{n=3}^{m} S_n^{},\varphi_{B,m}^{},\delta_B^{},\rho_{B,m}^{})$.

\section{Justification of the algorithms A and B}
\label{sectproof}

\begin{theorem}
$A_m$ is an ubr-algorithm for $\P_m$ if $m \geq 2$, $B_m$ is an ubr-algorithm
for $\P_m$ if $m \geq 3$.
\end{theorem}

To prove this, we have to verify the four conditions that we required in section
\ref{sectprin}. The second one is fulfilled by construction, because in the definition of $\delta$
we used reducing moves. There can be at most $\#S-1$ reducing steps for elements
in a linearly ordered, finite set $S$. So $j=\#S-1$ is an integer satisfying the
condition $\delta^j = \delta^{j+1}$.
The remaining parts of the proof of this theorem are given in the rest of this section.

\begin{corollary}
\label{corineq}
We have the inequalities
$$O_A(m)\;\geq\; \rk \P_m \;\geq\; O_C(m)  \;\;\;(\mbox{for }m \geq 2) \mbox{\quad and\quad }
O_B(m) \;\geq\; \rk \P_m \;\;\;(\mbox{for }m \geq 3)$$
where $O_A, O_B, O_C$ are given by $\;O_A(m)\; := \;\mbox{Output}(A_m)$, $\;O_B(m) \;:= \;\mbox{Output}(B_m)$ and
$O_C(m)\; := \;\dim \tilde \Phi_m \big(\Zmodule{\mbox{ caterpillar diagrams of degree } m\; }\big)$.
\end{corollary}

\subsection{Verification of the first condition}

In a first step we show that already the simply connected circle
diagrams generate ${\cal P}$.
One should recall that the circle on which the univalent vertices have been glued
is just a means of visualization, not a part of the diagram.
By a {\it loop} of a diagram we mean a closed path of consecutive edges
that meets each edge at most once.
Obviously, a loop contains only trivalent vertices and it cannot encounter a
vertex twice.
A vertex is called {\it bound} if there exists a loop going through that vertex,
otherwise it is called {\it free}.

We have a threefold partition of circle diagrams, given by their degree,
the dimension of the first homology and the number of free trivalent
vertices: $$D_{m,k,n} := \{\; D \;\vert\; D \mbox{ has }
2m \mbox{ vertices, } \dim H_1(D) = k,\; n \mbox{ of the trivalent vertices are free }\}$$

\begin{lemma}
\label{lemmatrees}
Any element of $D_{m,k,n}$ with $k>0$ and $n>0$ can be expressed (over $\P_m$) in terms
of elements of $D_{m,k,n-1}$.
\end{lemma}
\proof
Any $D \in D_{m,k,n}$ has by assumption both bound and free trivalent vertices.
$D$ is connected, so there exists an edge
connecting a free trivalent vertex $f$ with a bound vertex $b$.
Then the application of the IHX relation on this edge
presents $D$ as difference of two diagrams $\in D_{m,k,n-1}$:
\cepsfig{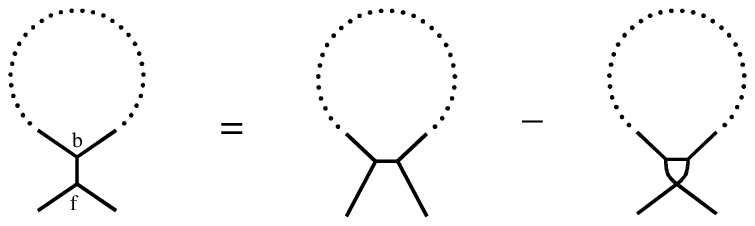}

\noindent
The IHX does not change the homology, but the diagrams on right side
have one free trivalent vertex less. \qed

\begin{lemma}
\label{lemmatrees2}
Any element of $D_{m,k,0}$ with $k>0$ can be expressed (over $\P_m$) in terms
of elements of $\bigcup_{i} D_{m,k-1,i}$.
\end{lemma}
\proof
Any circle diagram has at least one free vertex because, by definition, it
has one or more univalent vertices. For $k > 0$ any $D \in D_{m,k,0}$ has
bound vertices. There must be an edge connecting a free vertex $f$ with
a bound vertex $b$. By assumption $f$ has to be a univalent vertex.
An application of the STU-relation at $f$ opens the loops going through $b$
without introducing new loops:
\cepsfig{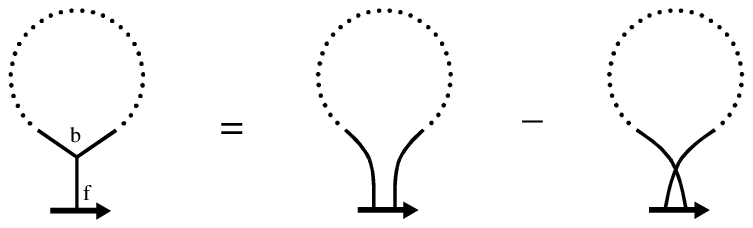}

\noindent
The diagrams $D_1, D_2$ on the right satisfy $\dim H_1 (D_i) = \dim H_1 (D)-1$. \qed

\begin{proposition}
\label{propphia}
$\{D^A_\pi\;\vert\;\pi\in S_{m-1}\}$ generates $\P_m$.
\end{proposition}
\proof
Lemmas \ref{lemmatrees} and \ref{lemmatrees2} imply that
the simply connected diagrams $D_{m,0,m-1}$ generate $\P_m$ 
(simply connected diagrams of degree $m$ must have $m-1$ trivalent vertices).

For $D \in D_{m,0,m-1}$ we choose two neighbouring (with respect to the cyclic ordering)
univalent vertices $a$, $b$. Let $P$ denote the (uniquely determined) path
in $D$ from $a$ to $b$.

If there are trivalent vertices that do not lie on
$P$, we can use exactly the same argument as in lemma \ref{lemmatrees} to
increase the number of vertices on $P$. We finally end up with diagrams
that have a path $P$ going through all trivalent vertices and that connects 
two neighbouring univalent vertices.
Because of the AS relation, we can sippose that all $m-1$ edges branching off $P$ 
lie on the left side of $P$. All these diagrams are of the
form $D^A_\pi$ with $\pi\in S_{m-1}$.\qed

\begin{proposition}
\label{propphib}
$\{D^B_\pi\;\vert\;\pi\in S_m\}$ generates $\P_m$ for $m \geq 3$.
\end{proposition}
\proof
We want to show that $D_{m,1,0}$ is a set of generators.
Because of lemma \ref{lemmatrees} and proposition \ref{propphia}, we only have to
show that every simply connected circle diagram can be expressed in terms of
diagrams having one loop, i.e.~elements of $\bigcup_i D_{m,1,i}$.
Let $D \in D_{m,0,m-1}$ with $m \geq 3$, then $D$ has a
trivalent vertex $t$ that is connected to two univalent vertices $a, b$.
We can "throw out" all other univalent vertices between $a$ and $b$ on the
circle, by using STU relations:
\cepsfig{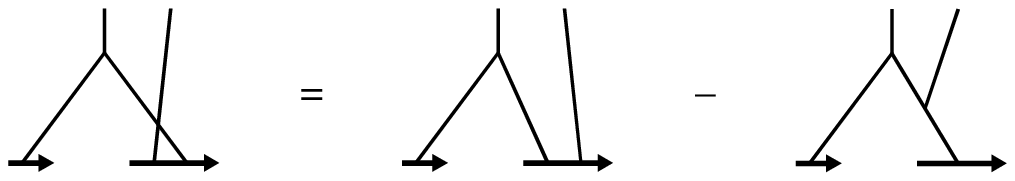}

\noindent
The second diagram on the right has a loop, so it remains to show that
a diagram with a trivalent vertex that is connected to two neighbouring
univalent vertices and another trivalent vertex, is equivalent to a diagram
with a loop. This is done by the following observation:
\cepsfig{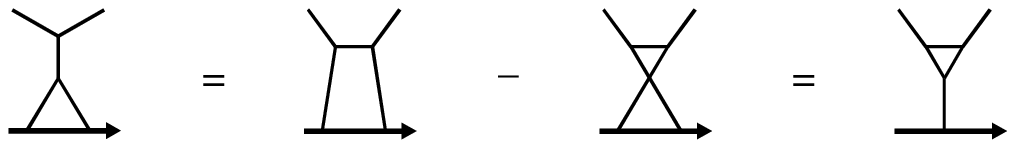}

\noindent
We have shown that diagrams having one loop and no free trivalent vertices
generate $\P_m$ for $m \geq 3$. These diagrams are equivalent
by the AS relation to $\pm D^B_\pi$ for some $\pi \in S_m$. \qed

\subsection{Verification of the third condition}

\begin{proposition}
The maps $\delta_A$ and $\delta_B$ induce the identity in $\P$.
\end{proposition}
\proof
We verify this for every move separately.
\begin{itemize}
\item[I:]
If the left and right side of this move are called $\pi$ and $\pi^\prime$,
then $D^A_\pi$ and $D^A_{\pi^\prime}$ have little triangles on the upper and lower end
(two of the sides are edges of the diagram, the third is part of the circle).
The second picture in the proof of proposition \ref{propphib}
together with the fact that $t \in \Lambda$, implies that
one can push the triangle down through the whole diagram, so $D^A_\pi = D^A_{\pi^\prime}$.
\item[II:] This move is an application of the STU relation; the triangle
in the second term on the right has been pushed down.
\cepsfig{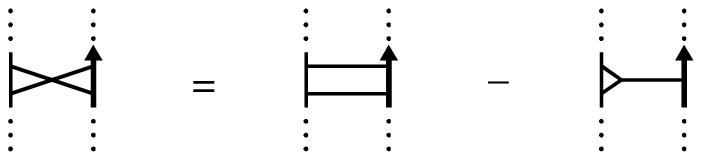}

\item[III:]
The four moves of type III occur, when the following diagram identities are resolved by
STU on the right and IHX on the left in two different ways that are
indicated by the small arrows.
\cepsfig{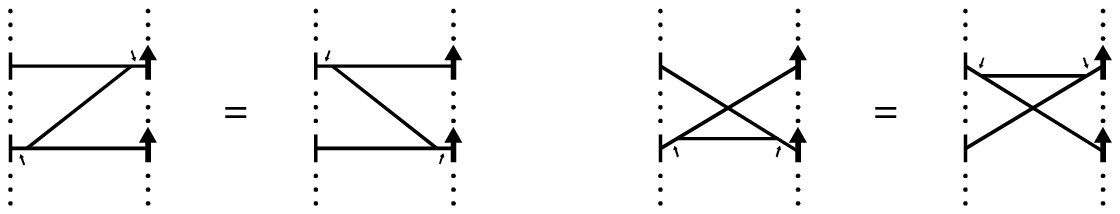}

\item[I$^\prime$:]
The diagrams $D^B_\pi$ are given by putting the permutation $\pi \in S_n$ between
a loop and an oriented circle. The multiplication with $\nu_n$ on the left or
right does not change the diagram, because it just slides vertices on
the loop or on the circle. The multiplication with $\mu_n$ corresponds to
a flip of the loop. When $n$ is odd the AS relation causes the sign to
change during the flip.
\item[II$^\prime$:]
The relation is of the form $\pi \rightarrow \pi_1 - \pi_2$ with
$\pi,\pi_1 \in S_n$ and $\pi_2 \in S_{n-1}$. The picture for type II
makes clear that $D^B_\pi = D^B_{\pi_1} - t\;D^B_{\pi_2}$.
\end{itemize}
Since $\delta_A$ and $\delta_B$ are defined by these moves, we have
hereby shown that $\hat\varphi_A\circ\delta_A = \hat\varphi_A$ and $\hat\varphi_{B,m}\circ\delta_B = \hat\varphi_{B,m}$.
\qed

\subsection{Verification of the forth condition}

First we have to understand $\Theta_k$.
For that purpose, we draw pictures with rectangular boxes named
$\Theta_k$ with $k$ entries on the left and $k$ exits on the right.
In each such box the pictures of all permutations occurring in
$\Theta_k$ shall be inserted (forgetting the upper $n-k$ constant strands)
 and the sum over all resulting diagrams (with the given signs) is taken.
In this way a picture with a $\Theta_k$-box in fact represents a linear
combination of $2^k$ diagrams.

\begin{lemma}
\label{lemmatheta}
In ${\cal A}$ the following relation holds:
\cepsfig{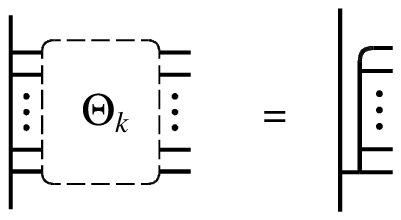}
\end{lemma}

\noindent
\proof (by induction on $k$)

\epsfxsize = 13cm  \epsfig{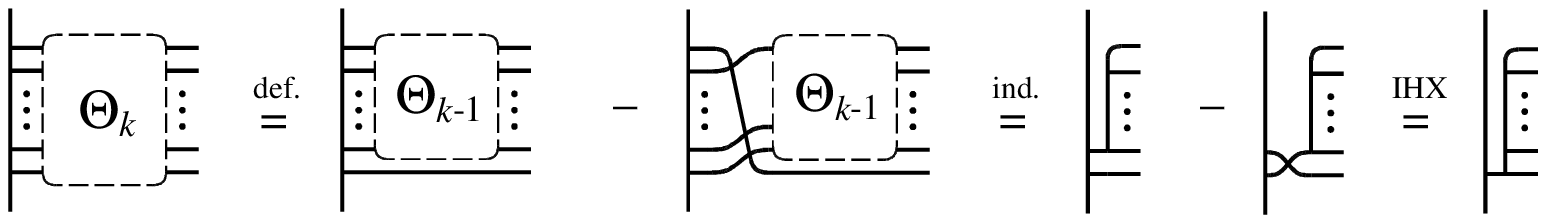} \qed

\begin{proposition}
$\hat\varphi_A(\pi) = (-1)^{n-k}\hat\varphi_A(\Theta_k \pi^\prime)$ for any
$\pi \in S_n$.
\end{proposition}
\proof
The external vertex on the lower side of $\hat\varphi_A(\pi)$ is named A, the
one on the upper side B and the lowest on the right side C. We will
rotate the circle clockwise, moving C $\rightarrow$ A $\rightarrow$ B.
This operation looks like this (in the picture of $\pi$
the line going from $\pi^{-1}(1)$ to $1$ has been drawn, the other lines have been omitted):
\cepsfig{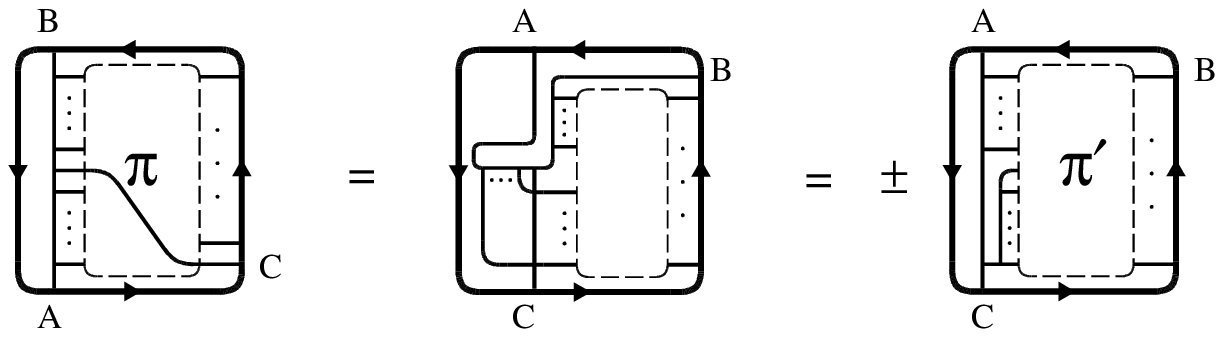}

\noindent
To get the second equality, one has to pull straight the path from C to A.
To do this the $n-k = \pi^{-1}(1)-1$ lowest trivalent vertices are swapped, which
is the reason for the factor $(-1)^{n-k}$.
The permutation in the box on the right side is $\pi^\prime$ of definition \ref{defrhoa}.
So by lemma \ref{lemmatheta} the third diagram in the equation is
equivalent to $\hat\varphi(\Theta_k \pi^\prime)$. \qed

\begin{lemma}
\label{lemmups}
For any $\pi\in S_n$, $0 \leq k \leq n-3$ the element $\hat\varphi_{B,n} \big(\Upsilon_{n,k}(\pi)\big)$ of $\P_{n+1}$
is equivalent to the following diagram.
\cepsfig{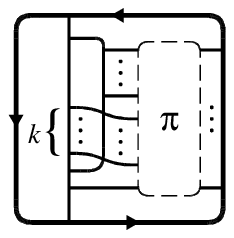}
\end{lemma}
\proof
The permutation $\pi^\# \in S_{n+1}$, which is used in the definition of $\Upsilon$
in section \ref{sectone}, is obtained
by doubling the $n-th$ string of $\pi$. So the two terms of $(1-\tau_n)\pi^\#$
allow a STU relation, after the right endpoints have been glued to the circle.
Together with our knowledge of $\Theta_k$, we get the following
diagram for $\chi_{1,n-k-2}^{} \; \Theta_{n-k-2}^{} \; \chi_{n-k-2,k+1}^{} \; (1 - \tau_n^{})\; \pi^\#$:
\cepsfig{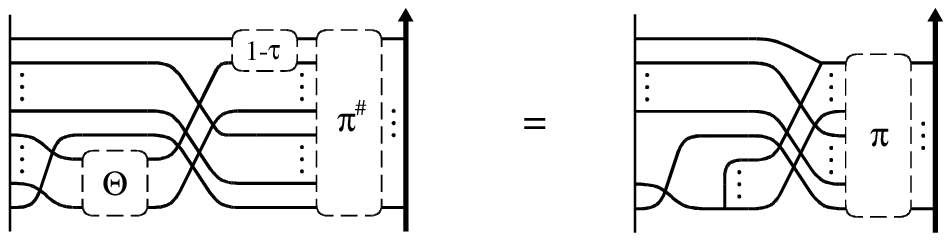}

\noindent
This results in the diagram of the claim.\qed

\begin{proposition}
$\hat\varphi_{B,m}^{}\big(\rho_{B,m}^{}(\pi)\big) = 0$ for all $\pi \in S_n$ with $3 \leq n \leq m$.
\end{proposition}
\proof
Let $x := (2\;1)(3\;5)(4) + (1)(2\;4)(3)$, then two moves of type I yield
$$\hat\varphi_{B,5}^{}(x) = -\hat\varphi_{B,5}^{}({\rm id}_5) - \hat\varphi_{B,5}({\rm id}_4)$$
Making three moves of type II$^\prime$ we get
$$\hat\varphi_{B,5}^{}(x) = \hat\varphi_{B,5}^{}\big((1)(2)(4)(3\;5)\big) = \hat\varphi_{B,5}^{}\big((1)(2)(3\;4\;5) - {\rm id}_4 + {\rm id}_3\big)$$
Both equations together imply
$\hat\varphi_{B,m}^{}\big((1)(2)(3\;4\;5)+{\rm id}_5-2{\rm id}_4+{\rm id}_3\big) = 0$
for $m \geq 5$.

\newln
For any $\pi \in S_m$ and $1 \leq q < m$, the STU relation allows
us to write $\hat\varphi(\pi) - \hat\varphi(\pi\tau_p)$ as a single diagram $D$ with
$\dim H_1(D) = 2$. For the second case in the definition of $\rho_{B,m}^{}$,
we keep the newly introduced loop of $D$ and eliminate the old one.
We show how this is done in the example $m = 9, p = 3, q = 7$:
\cepsfig{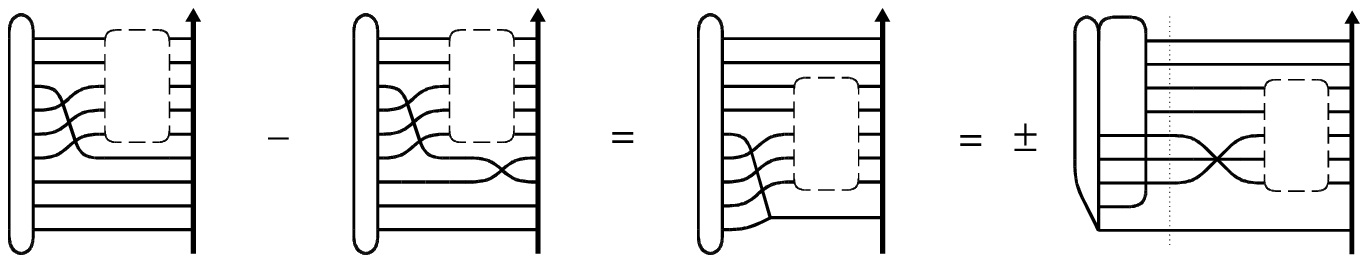}

\noindent
During this operation the diagram $D$ is not changed, except that the
$q-p-1$ trivalent vertices at between $p$ and $q$ are swapped with
the AS relation.
The permutation on right of the dotted line is $\pi^\prime \in S_{m-1}$
given in definition \ref{defrhob}. So by lemma \ref{lemmups}, the last diagram in the equation
is $\Upsilon_{n,q-p-1}(\pi^\prime)$, and we get
$\hat\varphi_{B,m}\big(\pi - \pi\tau_p +(-1)^{q-p}\Upsilon_{n,q-p-1}(\pi^\prime)\big) = 0$.

\newln
If $\pi \in S_n$ with $3 \leq n < m$ and $\pi (1)=1$, let $D_0 := D^B_\pi$.
Let the $D_1$ be the result of inserting a triangle at the lowest trivalent vertex of $D_0$.
Then $\varphi(\pi) = t^{m-n}D_0 = t^{m-n-1}D_1$.
Lemma \ref{lemmups} shows that $D_1 = \Upsilon_{n,0}(\pi)$, and so we have
$\hat\varphi_{B,m}\big(\pi - \Upsilon_{n,0}(\pi)\big) = 0$. \qed

\section{Remarks to the implementation}
\label{sectimpl}
\subsection{Dimensions}
One reason for the success of the presented algorithms is that the dimension of ${\cal I}$
($=$ number of irreducible permutations) is surprisingly small.
The moves we are using are very powerful, in the sense that
only a small number of permutation survive the reductions.
We tried out a large number of additional moves, but no considerable
improvement has been achieved this way.

\noindent
The following table displays 
%the number of connected circle diagrams (modulo the AS relation), 
the number of permutations and the number of irreducible
permutations, for the ubr-algorithms $A$ and $B$:

\smallskip
\noindent
{
\footnotesize
\begin{tabular}{||c|c|c|c|c|c|c|c|c|c|c|c|c||}
\hline
m                     &$ 2 $&$  3 $&$ 4 $&$  5 $&$   6 $&$   7 $&$    8 $&$     9 $&$     10 $&$      11 $&$       12 $\\
%diagrams              &$ 1 $&$  2 $&$   $&$    $&$     $&$     $&$      $&$       $&$        $&$         $&$          $\\
\hline
$\dim        S_A$     &$ 1 $&$  2 $&$ 6 $&$ 24 $&$ 120 $&$ 720 $&$ 5040 $&$ 40320 $&$ 362880 $&$ 3628800 $&$ 39916800 $\\
$\dim {\cal I}_A$    &$ 1 $&$  1 $&$ 2 $&$  5 $&$  16 $&$  64 $&$  301 $&$  1583 $&$   9145 $&$   57449 $&$   389668 $\\
\hline
$\dim        S_B$     &$   $&$  6 $&$ 30 $&$ 150 $&$ 870 $&$ 5910 $&$ 46230 $&$ 409110 $&$ 4037910 $&$ 43954710 $&$ 522956310 $\\
$\dim {\cal I}_B$    &$   $&$  1 $&$ 2 $&$  5 $&$  10 $&$  24 $&$  78 $&$  331 $&$  1685 $&$   9589 $&$   59782 $\\
\hline
\end{tabular}
}
\smallskip

\newln
The second reason for the success is that $\rho$ is "complicated enough"
to reproduce the kernel of $\hat\varphi$. It should not surprise that 
the "correct" $\rho_{A}^{}$ and $\rho_{B,m}^{}$ have been 
found by an intensive trial and error process.
We did not expect that the calculated upper bounds are sharp;
in fact, the algorithms described here are modifications of parts of a
much bigger program that computed the "exact" value $\rk \P$.

\subsection{Hints to the implementations}
\label{hintsec}

At first a list of irreducible permutations for the desired degree $m$ should
be made. Then $\rho(\pi)$ is calculated for any $\pi$ in this list.

The real difficulty is to evaluate $\bar\rho(\pi) = \Delta(\rho(\pi))$.
The simplest idea is to consecutively apply $\delta$, until a linear
combination of irreducible elements is reached. But this would be much too
slow for the interesting degrees $> 9$, because the
reduction trees are too nested.

The solution is to do it upside down. After assigning values to the minimal
permutations, we go from small permutations to the bigger ones. If
we know the values of all permutations smaller than $\pi$, then
the value of $\pi$ is given by a single application of $\delta$ and
picking at most $7$ values out of the table.

One has to assign a $\dim {\cal I}$-dimensional vector to every $\pi$.
A short look at the dimensions shows that keeping this table in memory exceeds
the capacity of any computer. But one can do the calculation component
per component.
Our implementation does $32$ components at a time,
yielding $32$ rows of the matrix for $\bar\rho$ in each run.

Even if each matrix entry uses only one bit, a file in which the matrix
for $\bar\rho_{B,12}^{}$ is saved contains 426 megabytes.
This is one reasons why we are working with $k = \F_{2}$.
The other reason is that, some time ago, we found diagrams 
$x \in \A_{10}$ for which we could show $2x = 0$ but not $x = 0$.
To find out, wether $A_{10}$ contains elements of order $2$ or not, was 
the main stimulus to make these computer computations.

The program that computes then rank and nullity of the matrices is a standard
Gaussian algorithm, which can of course be implemented very efficiently for
$k = \F_{2}$.
By the way, the matrices occurring are not at all sparse: about 40\% of the entries are $1$ and
gzip compresses the files only by factors about 0.9.

It is not a bad idea to add a check sum to each row/column in the data
files, because the probability of making an error in reading/writing a bit to
hard disc might (on some systems) not be far enough away from $1:10^{11}$, which is the approximate number
of bits that have to be read (in our implementation).

\section*{Acknowledgment}
I would like to thank C.-F.~B\"odigheimer, Michael Eisermann, Christoph Lamm
and Jens Lieberum for helpful discussions and reading the manuscript.
I am also thankful to the Studienstiftung for financial support and
the Graduiertenkolleg for providing the computer on
which the computation has been performed.

\section*{Computer files}

At the following locations in the internet,

{\tt http://www.uni-bonn.de/\~{}jk/pvi12.html}

{\tt ftp://ftp.uni-bonn.de/usr/jk/pvi12}

\noindent
you will find:
\begin{itemize}
\item
C/Pascal implementations of the algorithms.
\item
a basis for $\P_{m\leq 12}$ that is produced by the algorithm $B$.
\item
the values of $\tilde\Phi$ of the 175 diagrams that are 
mentioned in the proof of theorem \ref{theop}.
\item
a file containing this paper (or a newer version of it).
\item
a summary of the results of section \ref{sectres} and possibly some newer data.
\end{itemize}

\newln
\newln
\noindent
{\tt Mathematisches Institut}

\noindent
{\tt Beringstrasse 1}

\noindent
{\tt D-53115 Bonn}

\noindent
{\tt e-mail:jk@math.uni-bonn.de}

\noindent
{\tt http://rhein.iam.uni-bonn.de/\~{}jk}

\end{document}